\documentclass[twocolumn]{autart}    

 \usepackage{subfig}
\usepackage{mathrsfs}
\usepackage{amsfonts}
\usepackage{graphicx}          
 \usepackage{subfig}
 \usepackage{color}
 \usepackage{amssymb}

\newtheorem{lemma}{Lemma}
\newtheorem{definition}{Definition}
\newtheorem{theorem}{Theorem}
\newtheorem{remark}{Remark}

\newtheorem{example}{Example}

\newtheorem{corollary}{Corollary}                               

\begin{document}

\begin{frontmatter}

\title{Controllability of {Kronecker Product Networks}}

\thanks[footnoteinfo]{This work is supported by the National Natural Science Foundation (NNSF) of China under Grants 11802006, {11572015}, U1713223, 61673026, {the Beijing Natural Science Foundation under Grant 1194024, the Fundamental Research Funds for the Central Universities under Grant YWF-19-BJ-J-275,} and the Hong Kong Research Grants Council under the GRF Grant CityU 11200317. The material in this paper was not presented at any conference.}

\author[YUQING]{Yuqing Hao}\ead{haoyq@buaa.edu.cn},    
\author[YUQING]{Qingyun Wang}\ead{nmqingyun@163.com},
\author[Duan]{Zhisheng Duan}\ead{duanzs@pku.edu.cn},             
\author[CHEN]{Guanrong Chen}\ead{eegchen@cityu.edu.hk}  

\address[YUQING]{Department of Dynamics and Control, Beihang University, Beijing, 100191, China}  
\address[Duan]{State Key Laboratory for Turbulence and Complex Systems, Department of Mechanics and Engineering Science, College of Engineering, Peking University, Beijing 100871, China}             
\address[CHEN]{Department of Electronic Engineering, City University of Hong Kong, Hong Kong, China}        


\begin{abstract}                          
A necessary and sufficient condition is derived for the controllability of {Kronecker product networks}, where the factor networks are general directed graphs. The condition explicitly illustrates how the controllability of the factor networks affects the controllability of the composite network. For the special case where at least one factor network is diagonalizable, an easily-verifiable condition is explicitly expressed. Furthermore, the controllability of higher-dimensional multi-agent systems is revisited, revealing that some controllability criterion reported in the literature does not hold. Consequently, a modified necessary and sufficient condition is established. The effectiveness of the new conditions is demonstrated through several examples.
\end{abstract}

\end{frontmatter}

\section{Introduction}
Controllability is a fundamental issue to be addressed before considering how to control a dynamical system in applications \cite{chen_pinning_2016}. This subject has been extensively studied over more than half a century with various criteria developed, such as the PBH test, Kalman and other kinds of rank conditions, substantial graphic properties, and so on \cite{hautus_controllability_1969,kalman_canonical_1962,Trentelman_2012}.


For most large-scale networked systems, these criteria cannot be applied {practically} because of {their} complex structures and heavy computational burdens. Therefore, in recent years, the notion of network controllability has received compelling attention with some efficient criteria established {\cite{liu_controllability_2011,controllability_circulant_2013,Notarstefano,controllability_path_2012,controllability_multichain_2017,partition_controllability_2017,antagonistic_controllability_2017,xuemengran_comment_2018,wang_controllability_2016,hao_controllability_2017,xuemengran_input_2017,xuemengran_model_2018,yuan_controllability_2013,caoming_2014,zhangyuan_2017,zhoutong_2015}. In \cite{yuan_controllability_2013}, an exact controllability framework was introduced to identify the minimum set of input nodes for a general network with an arbitrary link-weight distribution. The controllability of networks with specific topologies, such as path graphs, cycle graphs, circulant graphs, multi-chain graphs and grid graphs, was explored in \cite{controllability_circulant_2013,Notarstefano,controllability_path_2012,controllability_multichain_2017}. Network controllability was investigated from a graph-theoretic perspective in \cite{Rahmani_2009,caoming_2014,partition_controllability_2017,antagonistic_controllability_2017}. It is noted that most of the above results are derived for the networks with one-dimensional nodes. Recently, the controllability of networks with
 higher-dimensional nodes has attracted a great deal of interest \cite{hao_controllability_2017,wang_controllability_2016,xuemengran_comment_2018,xuemengran_input_2017,xuemengran_model_2018,caoming_2014,zhangyuan_2017,zhoutong_2015}. Some controllability conditions for networked LTI systems were presented in \cite{zhoutong_2015,zhangyuan_2017}, which depend on the transmission zeros of every subsystem and the connection matrix. The controllability of diffusively coupled LTI systems was studied in \cite{caoming_2014} and \cite{xuemengran_comment_2018}. In \cite{wang_controllability_2016}, the controllability condition for networked MIMO systems was established in terms of two algebraic matrix equations. Moreover, some easily-verifiable conditions were proposed in \cite{hao_controllability_2017,xuemengran_input_2017} for the controllability of networked LTI systems with a diagonalizable topology matrix.}

Besides isolated networks, composite networks have come {into} play due to their broad applications in {different areas of engineering \cite{book,Loan_2000}}. There are many kinds of composite networks, such as Cartesian product networks, Kronecker product networks, strong product networks, {lexicographic product networks} and so on. {In addition to their importance of constructing `larger' networks out of `small' ones, they are useful in the sense that one can get insights about the properties of composite networks from the factor networks.} Intuitively, the controllability of {a} composite network might be verified by checking some properties of the factor networks. The controllability and observability of Cartesian product networks were investigated in \cite{chapman_controllability_2014}. In \cite{Notarstefano}, the controllability and observability of linear dynamical systems whose dynamics are induced by the Laplacian of a grid graph were studied. {Note that many real-world networks are similar to stochastic Kronecker product graphs. For example, the Kronecker power of a simple generating matrix can yield a Kronecker product graph that fits the Internet (at the autonomous systems level) fairly well \cite{leskovec_2010}. Large online social networks, web and blog graphs, peer-to-peer networks, etc., can also be modeled by Kronecker product networks. Moreover, every non-trivial graph has a prime factorization over the Kronecker product \cite{book,Loan_2000}. It has lower computational complexity to check the controllability of a large-scale network by examining some properties of the smaller factor networks. Therefore,} the controllability analysis for Kronecker product networks has brought about renewed interest recently \cite{Asavathiratham_2001,chapmanconference,xuemengran,horn_matrix_1987}. {The eigenanalysis for the Kronecker product of two matrices was presented in \cite{horn_matrix_1987}, which shows that the Kronecker products of the factor matrices' eigenvectors are the eigenvectors of the composite matrix. However, not all the eigenvectors of the composite matrix are formulated therein.} In \cite{xuemengran}, the Kronecker product of defective matrices was revisited, characterizing the number of the eigenvectors. {Nevertheless,} the explicit expressions of the eigenvectors, which are the cornerstone of controllability analysis, are not presented. Recently, a sufficient condition was established in \cite{chapmanconference} for the controllability of Kronecker product networks, where {the topology matrices are} required to be diagonalizable.

In this paper, the controllability of {Kronecker product networks} is revisited. The contribution of the paper is four-fold. First, the factor networks considered here {are} general, directed and weighted. Differing from the condition given in \cite{chapmanconference}, which requires the {topology} matrix of the composite network to be diagonalizable, this paper removes the diagonalizability requirement. Second, a new necessary and sufficient condition on the controllability of Kronecker product {networks} is provided in terms of eigenvectors. Compared with the classical PBH test, the new condition typically has a much lower computational cost. Third, for the special case where at least one factor network is diagonalizable, a specified condition is explicitly expressed, which is easier to verify. Finally, this paper shows that the sufficiency of the controllability criterion given in \cite{caining} does not hold, thereby a modified necessary and sufficient condition is derived.

The remainder of this paper is organized as follows.
Some notations and preliminaries are given
in Section \ref{section2}. The model is formulated in Section \ref{section3}. Some conditions on the controllability of {Kronecker product networks} are developed in Section \ref{section4}. The controllability of higher-dimensional multi-agent systems is reinvestigated in Section \ref{section5}. Finally, conclusions are drawn in Section \ref{section6}.

\section{Notations and Preliminaries}
\label{section2}
In this section, notations and useful preliminaries are introduced.

\subsection{Notations}
{The notations mostly follow \cite{hao_controllability_2017}. Let $e_i$ be the row vector with all zero entries except for $[e_i]_i=1$.} The linear span of row vectors $v_1$, $\cdots$, $v_k$ is {a set of all the linear combinations} of these vectors, i.e., $span\{ {v_1},{v_2}, \cdots ,{v_k}\}  = \left\{ {\sum\limits_{i = 1}^k {{c_i}{v_i}} {\rm{|}}{{\rm{c}}_i} \in \mathbb{R}} \right\}$. Let ${U_1} \oplus {U_2}$ be the direct sum of two spaces $U_1$ and $U_2$.
Matrices, if their dimensions are not explicitly indicated, are assumed to be compatible for algebraic operations.

{A weighted digraph ${\mathcal{G}}=(V,E,W)$ is characterized by a node set $V$ with cardinality $n$, an edge set $E$ comprised of ordered pairs of nodes with cardinality $m$, and a weight set $W$
with cardinality $m$, where an edge exists from nodes $i$ to $j$ if
$(i,j)\in E$ with edge weight $w_{ji}\in W$.
${\mathcal A}=[a_{ij}] \in \mathbb{R}^{n\times n}$ is called the adjacency matrix of $\mathcal{G}$ with $a_{ij}=w_{ij}$ if $(j,i) \in E$ and $a_{ij}=0$ otherwise. The $i$th diagonal term of the adjacency matrix denotes the weight of the self-loop on node $i$. }

\subsection{{Kronecker Product Network}}
\begin{figure}[htb]
\centering
\subfloat[${\mathcal G}_1$]{%
\includegraphics[width=.15\textwidth]{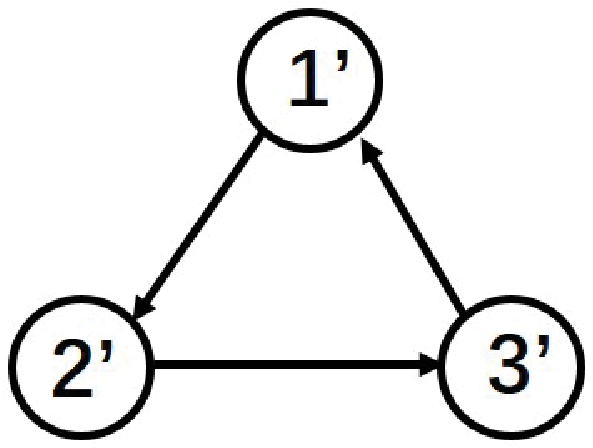}}
\subfloat[${\mathcal G}_2$]{%
\includegraphics[width=.15\textwidth]{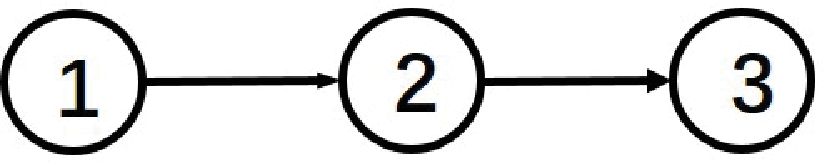}}\hfill
\subfloat[${\mathcal G}={\mathcal G}_1 \times {\mathcal G}_2$]{%
\includegraphics[width=.15\textwidth]{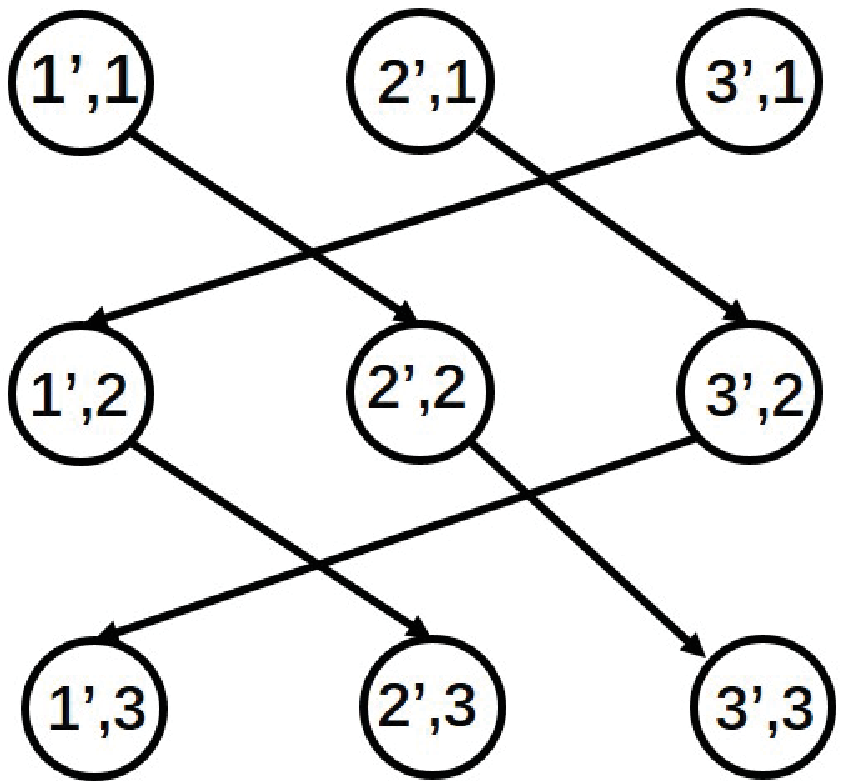}}\\
\caption{{Factor graphs ${\mathcal G}_1$ and ${\mathcal G}_2$ and their Kronecker product ${\mathcal G}={\mathcal G}_1 \times {\mathcal G}_2$}}
\label{fig1}
\end{figure}

The Kronecker product network is a kind of composite {network}
that can be obtained by applying Kronecker product operation(s) to several smaller networks, called factor networks. Let $\mathcal{G}_1=(V_1,E_1,W_1)$ and $\mathcal{G}_2=(V_2,E_2,W_2)$ be two factor networks. The Kronecker product of $\mathcal{G}_1$ and $\mathcal{G}_2$, denoted by $\mathcal{G}=\mathcal{G}_1 {\times} \mathcal{G}_2$, has the node set $V_1 \times V_2$. There is an edge from
node $(i,p)$ to node $(j,q)$ if and only if $(i,j)$ is an
edge of $E_1$ and $(p,q)$ is an edge of $E_2$. If an edge exists, the corresponding
weight is $w_{((i,p),(j,q))} = w_{ij}w_{pq}$. An example of the Kronecker product of two graphs
$\mathcal{G}_1$ and $\mathcal{G}_2$ is displayed in Fig. \ref{fig1}.

\subsection{Useful Lemmas and Definitions}
%

\begin{lemma}\cite{roman_advanced_2005}
\label{lemma2}
If $u_1$, $u_2$, $\cdots$, $u_n$ are linearly independent vectors and $v_1$, $v_2$, $\cdots$, $v_n$ are arbitrary vectors, then
$\sum\limits_{k = 1}^n {{u_k} \otimes {v_k} = 0} $
implies that $v_k =0$ for all $k=1,2,\cdots,n$. Moreover, the roles of $u_k$ and $v_k$ in the above statement can be exchanged.

\end{lemma}

\begin{definition}\cite{roman_advanced_2005}
A row vector $x_m$ is called {the} $m$th-order generalized left eigenvector of matrix $A$ corresponding to its eigenvalue $\lambda$ if
${x_m}{(A - \lambda I)^m} = 0$,
and
${x_m}{(A - \lambda I)^{m - 1}} \ne 0$.
Moreover, $x_1$, $\cdots$, $x_g$ form a left Jordan chain of $A$ on top of $x_1$, where the maximum value of $g$ is called the length of this Jordan chain.

\end{definition}

\section{Model Formulation}
\label{section3}

{Consider a network consisting of $Nn$ nodes with a directed and weighted topology {${\mathcal G}$} in the following form \cite{chen_book,controllability_metrics_2014,leskovec_2010}:
\begin{equation}
\label{eq1}
 \begin{array}{l}
 {{\dot x}_i}(t) =  \sum\limits_{j = 1}^{Nn} {{c _{ij}}{x_j}(t) + {\delta_i}{u_i}(t),}\;\;i = 1,2, \cdots ,Nn, \\
 \end{array}
\end{equation}
 where ${x_i} \in {\mathbb{R}}$ is the state of node $i$, ${c _{ij}} \in \mathbb{R}$ represents the coupling strength between nodes $i$ and $j$, ${u_i} \in {\mathbb{R}}$ is the control input to node $i$, and ${\delta_i} = 1$ if node $i$ is under control, but otherwise ${\delta_i}=0$, for all $i = 1,2, \cdots ,Nn$. Assume that ${c _{ij}} \ne 0$ if there is an edge from node $j$ to node $i$, otherwise ${c _{ij}} = 0$, for all $i,j = 1,2, \cdots ,Nn$.
Denote
$A({\mathcal G}) = {\mathcal A}({\mathcal G})= [{c_{ij}}] \in {\mathbb{R}^{{Nn} \times {Nn}}}$ and $B = diag\left\{{\delta_1},{\delta_2}, \cdots, {\delta_{Nn}}\right\}$,
which represent the topology and the external input channels of the network (\ref{eq1}), respectively.
Let $X = {[ {x_{\rm{1}},\;x_2,\; \cdots ,\;x_{Nn}} ]^T}$ be the whole state of the network, and $U = {[ {u_{\rm{1}},\;u_2,\; \cdots ,\;u_{Nn}} ]^T}$ be the total external control input. Then, the network (\ref{eq1}) can be rewritten in a compact form as
\begin{equation}
\label{eq3}
\dot X(t) = A({\mathcal G}) X(t) + B U(t).
\end{equation}
In the following, consider the controllability of the network with topology graph ${\mathcal G}$ being the Kronecker product of two factor graphs ${\mathcal G}_1$ and ${\mathcal G}_2$. The dynamics of the factor networks are described by
${{\dot Y}_1}(t) = A({{\mathcal G}_1}){Y_1}(t) + {B_1}{U_1}(t)$ and
${{\dot Y}_2}(t) = A({{\mathcal G}_2}){Y_2}(t) + {B_2}{U_2}(t)$,
where $Y_1 \in \mathbb{R}^{N}$ and $U_1 \in \mathbb{R}^{D}$ are the state vector and the control input of the first factor network, respectively; $Y_2 \in \mathbb{R}^{n}$ and $U_2 \in \mathbb{R}^{d}$ are the state vector and the control input of the second factor network, respectively. {According to the definition of Kronecker product graph, it is easy to prove $A({\mathcal G}) {=} A({{\mathcal G}_1} {\times} {{\mathcal G}_2})= A({{\mathcal G}_1}) {\otimes} A({{\mathcal G}_2})$}. Therefore,
for this Kronecker product network, its compact form can be formulated as (\ref{eq3}) with
\begin{equation}\label{eq4}
\begin{array}{l}
A({\mathcal G}) {=} A({{\mathcal G}_1} {\times} {{\mathcal G}_2})= A({{\mathcal G}_1}) {\otimes} A({{\mathcal G}_2}),\\
B = {B_1} \otimes {B_2}.
\end{array}
\end{equation}
The analysis here is presented in terms of two factor networks, which can be extended to larger numbers of factor networks by sequential compositions.
In the following,
it will be shown that the controllability of the composite network can be revealed by
examining some features of the smaller factor networks, which makes the computational complexity much lower.
\begin{remark}
A graph $\mathcal G$ is prime if it is nontrivial and cannot be decomposed into the Kronecker product of two nontrivial graphs. An expression ${\mathcal G}={\mathcal G}_1 \times {\mathcal G}_2 \times \cdots \times {\mathcal G}_k$, with each ${\mathcal G}_i$ being prime, is called a prime factorization of ${\mathcal G}$. Note that every nontrivial graph has a prime factorization over the Kronecker product \cite{book}. It has lower computational complexity to check the controllability of a large-scale network by examining some properties of the smaller factor networks. Therefore, exploring low-dimensional controllability conditions for Kronecker product networks is of great significance in engineering applications, which gives insight into the controllability of large-scale networks.
\end{remark}}

\section{Main Results}
\label{section4}
In this section, controllability conditions for the composite network (\ref{eq3})-(\ref{eq4}) are considered. Firstly, the general case that both $A(\mathcal{G}_1)$ and $A(\mathcal{G}_2)$ are non-diagonalizable is investigated. Then, a special case where at least one factor network is diagonalizable is further analyzed, with a direct and easily-verifiable condition derived.

\subsection{Both $A(\mathcal{G}_1)$ and $A(\mathcal{G}_2)$ are non-diagonalizable}
In this subsection, the general case that both $A(\mathcal{G}_1)$ and $A(\mathcal{G}_2)$ are non-diagonalizable is investigated. The connection between the eigenspaces of the factor networks and that of the composite network provides a mechanism to establish efficient controllability conditions.
Firstly, {all} left eigenvectors of the Kronecker product of two defective matrices are characterized.

\begin{theorem}
\label{theorem5}
Let $P=\lambda I_{p}+N_{p} \in \mathbb{C}^{p \times p}$ and $Q=\mu I_{q}+N_{q} \in \mathbb{C}^{q \times q}$ be two defective matrices, where $N_{p}$ and $N_{q}$ are nilpotent matrices. The eigenvalues of $P\otimes Q$ are $\lambda \mu$. Moreover,

\begin{itemize}
  \item if $\lambda \mu \ne 0$, the corresponding left eigenvectors are $\eta^1=e_{p} \otimes \xi_1$, $\eta^2=e_{p-1} \otimes \xi_1+e_{p} \otimes \xi_2$, $\cdots$, $\eta^{\theta}=e_{p-\theta+1} \otimes \xi_1+e_{p-\theta+2} \otimes \xi_2+\cdots+e_{p} \otimes \xi_{\theta}$, where $\xi_1=e_{q}$, ${\xi _k} = \frac{{{{( - 1)}^{k + 1}}}}{{\lambda^{k - 1}}}\sum\limits_{l = 0}^{k - 2} {C_{k - 2}^l\mu^{k - l - 1}e_{{q} - k + l + 1}} $, $k=2,\cdots,\theta$, $\theta  = \min \{ {p},{q}\} $;
  \item if $\lambda=0$ and $\mu \ne 0$, the corresponding left eigenvectors are $\eta^k=e_{p} \otimes e_k$, $k=1,\cdots,{q}$;
  \item if $\lambda \ne 0$ and $\mu = 0$, the corresponding left eigenvectors are $\eta^{k}=e_{k} \otimes e_{q}$, $k=1,\cdots,p$;
  \item if $\lambda = \mu = 0$, the corresponding left eigenvectors are $\eta^1=e_{p} \otimes e_1$, $ \eta^2=e_{p} \otimes e_2$, $\cdots$, $\eta^{q}=e_{p} \otimes e_{q}$, $\eta^{q+1}=e_{p-1} \otimes e_{q}$, $\cdots$, $\eta^{q+p-1}=e_{1} \otimes e_{q}$.
\end{itemize}

\end{theorem}

{\textbf{Proof:}}
According to the structure of matrix $P \otimes Q$,
it is easy to verify that the eigenvalues of $P \otimes Q$ are $\lambda \mu$. In the following, the aforementioned four cases are proved respectively.

Case $1$: $\lambda \mu \ne 0$. Since $\eta^1({P} \otimes {Q}) = (e_{{p}} \otimes e_{{q}})({P} \otimes {Q}) = ({\lambda}e_{{p}}) \otimes ({\mu }e_{{q}}) = {\lambda }{\mu }\eta^1$, it follows that $\eta^1$ is the left eigenvector of $P \otimes Q$.
Note that $\eta^2$ is the left eigenvector of $P \otimes Q$ if and only if
\begin{equation}
\label{eq11}
\eta^2({P} \otimes {Q}) = {\lambda}{\mu }\eta^2.
\end{equation}
This implies that ${\lambda}{\xi _2}({\mu}I - {Q}) = {\xi _1}{Q}$. Since ${\lambda }{\xi _2}({\mu }I - {Q}) = {\mu }e_{{q}}$ and ${\xi _1}{Q} = {\mu }e_{{q}}$, equation (\ref{eq11}) holds, indicating that $\eta^2$ is the left eigenvector of $P \otimes Q$.

Assume that $\eta^1$, $\cdots$, $\eta^{k-1}$ are the left eigenvectors of $P\otimes Q$. Then, $\eta^k$ is the left eigenvector of $P\otimes Q$ if and only if
\begin{equation}
\label{eq12}
{\lambda }{\xi _k}({\mu }I - {Q}) = {\xi _{k-1}}{Q}.
\end{equation}
Since ${\lambda }{\xi _k}({\mu }I - {Q}) = \frac{{{{( - 1)}^{k + 2}}}}{{\lambda^{k - 2}}}\sum\limits_{l = 0}^{k - 2} {C_{k - 2}^l\mu^{k - l - 1}e_{{q} - k + l + 2}} $ and
${\xi _{k - 1}}{Q} = \frac{{{{( - 1)}^k}}}{{\lambda^{k - 2}}}\sum\limits_{l = 0}^{k - 3}$ $C_{k {-} 3}^l\mu^{k {-} l {-} 2}({\mu}e_{{q} {-} k {+} l {+} 2} {+}  e_{{q} {-} k {+} l {+} 3})$
$ = \frac{{{{( - 1)}^k}}}{{\lambda^{k - 2}}}\sum\limits_{l = 0}^{k - 2} C_{k - 2}^l\mu^{k - l - 1}e_{{q} - k + l + 2},$
one can verify that (\ref{eq12}) holds. Therefore, $\eta^k$ is the left eigenvector of $P\otimes Q$.
It has been shown in \cite{xuemengran} that the number of the left eigenvectors for this case is $\min\{p,q\}$, thus all the left eigenvectors of $P\otimes Q$ are $\eta^1$, $\cdots$, $\eta^{\theta}$.

Case $2$: $\lambda=0$ and $\mu \ne 0$. Since $\eta^k({P} \otimes {Q}) = (e_{{p}} \otimes e_k)({P} \otimes {Q}) = (e_{{p}}{P}) \otimes (e_k{Q}) = \textbf{0}_{{p}}^T \otimes (e_k{Q}) = \textbf{0}_{{p}{q}}^T = {\lambda }{\mu }\eta^k$, it follows that $\eta^k=e_{p} \otimes e_k$ is the left eigenvector of $P \otimes Q$, $k=1,\cdots,{q}$.
It has been shown in \cite{xuemengran} that the number of the left eigenvectors for this case is $q$, thus all the left eigenvectors of $P\otimes Q$ are $\eta^1$, $\cdots$, $\eta^{q}$.

One can prove the results {in} Cases $3$ and $4$ similarly, thus the detail is omitted.
This completes the proof.

\hfill $\blacksquare$

\begin{remark}
The number of the eigenvectors of $P \otimes Q$ associated with the sole eigenvalue $\lambda\mu$ {was given in \cite{xuemengran}}. However, it does not present explicit expressions of the eigenvectors, which are characterized in the above theorem. Theorem \ref{theorem5} is the cornerstone of the following controllability analysis for Kronecker product {networks}.

\end{remark}

In what follows, the left eigenvectors of $A(\mathcal{G})$ are expressed through the generalized eigenspaces of the factor networks.

\begin{theorem}
\label{theorem6}
Let $\lambda_1,\lambda_2,\cdots,\lambda_s$ be the eigenvalues of $A(\mathcal{G}_1)$, and $\mu_1,\cdots,\mu_t$ be the eigenvalues of $A(\mathcal{G}_2)$. Then, the eigenvalues of $A(\mathcal{G})$ are $\lambda_1\mu_1$, $\cdots$, $\lambda_1\mu_t$, $\cdots$, $\lambda_s\mu_1$, $\cdots$, $\lambda_s\mu_t$. Moreover, for the eigenvalue $\lambda_i \mu_j$ with geometric multiplicity $\theta_{ij}$,
\begin{itemize}
  \item if $\lambda_i\mu_j \ne0$, {then} $\theta_{ij}=\min\{p_i,q_j\}$, {and} the corresponding left eigenvectors are $\eta _{ij}^k {=} v_i^k {\otimes} w_j^1 {-} \frac{{\mu _j}}{{{\lambda _i}}}v_i^{k {-} 1} {\otimes} w_j^2 \\{+} v_i^{k {-} 2} \otimes \left[ {\frac{{\mu _j^2}}{{\lambda _i^2}}w_j^3 {+} {\frac{{{\mu _j}}}{{\lambda _i^2}}}w_j^2} \right]+\cdots+ v_i^1 \otimes \left[ {\frac{{{{( {-} 1)}^{k {+} 1}}}}{{\lambda _i^{k {+} 1}}}{\sum\limits_{l {=} 0}^{{k {{-}} 2}}} {C_{{k{{-}}2}}^l\mu _j^{k{-}l{-} 1}}}\right.\\\left.{{w_j^{k {-} l}} } \right]$, $k=1,\cdots,\theta_{ij}$;
  \item  if $\lambda_i=0$ and $\mu_j \ne 0$, {then} $\theta_{ij}=q_j$, {and} the corresponding left eigenvectors are $\eta_{ij}^k=v_i^1 \otimes w_j^k$, $k=1,\cdots,\theta_{ij}$;
  \item if $\lambda_i\ne0$ and $\mu_j = 0$, {then} $\theta_{ij}=p_i$, {and} the corresponding left eigenvectors are $\eta_{ij}^k=v_i^k \otimes w_j^1$, $k=1,\cdots,\theta_{ij}$;
  \item if $\lambda_i=0$ and $\mu_j = 0$, {then} $\theta_{ij}=q_j+p_i-1$, {and} the corresponding left eigenvectors are $\eta_{ij}^1=v_i^1 \otimes w_j^1$, $\cdots$, $\eta_{ij}^{q_j}=v_i^1\otimes w_j^{q_j}$, $\eta_{ij}^{q_j+1}=v_i^2\otimes w_j^{1}$, $\cdots$, $\eta_{ij}^{\theta_{ij}}=v_i^{p_i}\otimes w_j^{1}$,
\end{itemize}
where $v_i^1$, $\cdots$, $v_i^{p_i}$ is the left Jordan chain of $A(\mathcal{G}_1)$ associated with the eigenvalue $\lambda_i$, and $w_j^1$, $\cdots$, $w_j^{q_j}$ is the left Jordan chain of $A(\mathcal{G}_2)$ associated with the eigenvalue $\mu_j$, $i=1,\cdots,s$, $j=1,\cdots,t$.

\end{theorem}

{\textbf{Proof:}}
Let $V \in \mathbb{C}^{N \times N}$ be a nonsingular matrix such that $VA(\mathcal{G}_1)V^{-1}{=}P{=}blockdiag\left\{P_1,\;P_2,{\cdots},P_s\right\}$, where $P$ is the Jordan form of $A(\mathcal{G}_1)$ and the $i$th Jordan block $P_i=\lambda_i I_{p_i}+N_{p_i}\in \mathbb{C}^{p_i \times p_i}$ with $N_{p_i}$ being a nilpotent matrix, $i=1,\cdots,s$. The left Jordan chain of $A(\mathcal{G}_1)$ associated with the eigenvalue $\lambda_i$ is denoted as $v_i^1$, $v_i^2$, $\cdots$, $v_i^{p_i}$, where $v_i^1$ is the top vector and $p_i$ is the length of the Jordan chain.

Let $W \in \mathbb{C}^{n \times n}$ be a nonsingular matrix such that $WA(\mathcal{G}_2)W^{-1}{=}Q{=}blockdiag\left\{Q_1,\;Q_2,{\cdots},Q_t\right\}$, where $Q$ is the Jordan form of $A(\mathcal{G}_2)$ and the $j$th Jordan block $Q_j=\mu_j I_{q_j}+N_{q_j}\in \mathbb{C}^{q_j \times q_j}$, $j=1,\cdots,t$. The left Jordan chain of $A(\mathcal{G}_2)$ associated with the eigenvalue $\mu_j$ is denoted as $w_j^1$, $w_j^2$, $\cdots$, $w_j^{q_j}$, where $w_j^1$ is the top vector and $q_j$ is the length of the Jordan chain.

Since $(V \otimes W)[A(\mathcal{G}_1) \otimes A(\mathcal{G}_2)]{(V \otimes W)^{ - 1}} = P \otimes Q = blockdiag\{ {P_1}\otimes{Q_1}, \cdots ,{P_1}\otimes{Q_t},$ $\cdots ,{P_s}\otimes{Q_1}, \cdots ,{P_s}\otimes{Q_t}\} $, the eigenvalues of $A(\mathcal{G})$ are $\lambda_1\mu_1$, $\cdots$, $\lambda_1\mu_t$, $\cdots$, $\lambda_s\mu_1$, $\cdots$, $\lambda_s\mu_t$.
Let ${V_i} = {[ {\begin{array}{*{20}{c}}
{v_i^{{p_i}T}}& \cdots &{v_i^{1T}}
\end{array}} ]^T}$ and ${W_j} = {[ {\begin{array}{*{20}{c}}
{w_j^{{q_j}T}}& \cdots &{w_j^{1T}}
\end{array}} ]^T}$. It is easy to derive that $({V_i} \otimes {W_j})A(\mathcal{G}) = ({P_i} \otimes {Q_j})({V_i} \otimes {W_j})$. If $\zeta   \in {\mathbb{C}^{1 \times {p_i} {q_j}}}$ is a left eigenvector of $P_i \otimes Q_j$, then $\zeta ({V_i} \otimes {W_j})$ is the left eigenvector of $A(\mathcal{G})$, $i=1,\cdots,s$, $j=1,\cdots,t$. Thus, the results follow from Theorem \ref{theorem5} directly.
\hfill $\blacksquare$

The explicit expressions of the left eigenvectors given in Theorem \ref{theorem6} are the core of the following controllability criteria for Kronecker product networks.
Let $U_{ij}=span\{\eta_{ij}^1,\cdots,\eta_{ij}^{\theta_{ij}}\}$ be the left eigenspace of $A(\mathcal{G})$ corresponding to the eigenvalue $\lambda_i \mu_j$.
In what follows, a theorem is established for the controllability of the Kronecker product network (\ref{eq3})-(\ref{eq4}).

\begin{theorem}
\label{theorem7}
Let $\Lambda=\{\lambda_1,\;\lambda_2,\;\cdots,\;\lambda_s\}$ be the set of the eigenvalues of $A(\mathcal{G}_1)$, and $\Sigma=\{\mu_1,\cdots,\mu_t\}$ be the set of the eigenvalues of $A(\mathcal{G}_2)$. The composite network (\ref{eq3})-(\ref{eq4}) is controllable if and only if the following two conditions hold simultaneously:
\begin{enumerate}
  \item $\forall \eta  \in {U_{ij}}$ and $\eta \ne 0$, $\eta (B_1  \otimes B_2) \ne 0$, for $i =1,\cdots,s$, $j =1,\cdots,t$;
  \item if $\lambda_{i_1}\mu_{j_1}=\lambda_{i_2}\mu_{j_2}=\cdots=\lambda_{i_r}\mu_{j_r}$, where $\lambda_{i_k} \in \Lambda$, $\mu_{j_k} \in \Sigma$, $k=1,\cdots,r$, $r >1$, then $\forall \eta  \in \mathop  \oplus \limits_{k = 1}^r {U_{{i_k}{j_k}}}$ and $\eta \ne 0$, $\eta(B_1  \otimes B_2) \ne 0$.
\end{enumerate}

\end{theorem}

{\textbf{Proof:}}
Using Theorem \ref{theorem6} and {the PBH test}, one can prove this theorem easily. Thus, the detail is omitted.
\hfill $\blacksquare$

\begin{remark}
Theorem \ref{theorem7} provides {a} precise and efficient {criterion} for determining the controllability of a general Kronecker product network, using the generalized eigenvectors of the factor networks with low dimensions. Compared with the classical PBH test, the new condition typically has a much lower computational cost. Due to linear systems duality, the results in Theorem \ref{theorem7} are equally applicable to the observability of Kronecker product {networks}.
\end{remark}

The effectiveness of the above condition can be illustrated by the following example.

\begin{example}
Consider graph ${\mathcal G}_1$ depicted in Fig. \ref{fig3}.
\begin{figure}[!htb]
  \centering
  \includegraphics[scale=0.5]{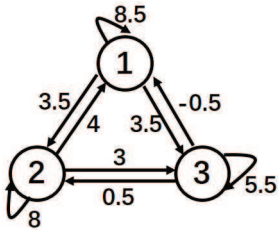}
  \caption{Graph ${\mathcal G}_1$}
  \label{fig3}
\end{figure}
It is easy to verify that
$A({{\mathcal G}_1}) = \left[ {\begin{array}{*{20}{c}}
   8.5 & 4 & -0.5  \\
   3.5 & 8 & 0.5  \\
   3.5 & 3 & 5.5  \\
\end{array}} \right]$.
The eigenvalues of $A({\mathcal G}_1)$ are $\lambda_1=12$, $\lambda_2=5$. The left eigenvector associated with $\lambda_1$ is $v_1=-e_1-e_2$ and the left Jordan chain corresponding to $\lambda_2$ is $v_2^1=e_2-e_3$, $v_2^2=-e_1+e_2$.
It is easy to obtain that $A({\mathcal G}_1)$ is non-diagonalizable.
Assume that the second node of ${{\mathcal G}_1}$ is under control.
Thus, one has
${B_1} = e_2^T$.
Let $A({\mathcal G})=A({\mathcal G}_1) \otimes A({\mathcal G}_1)$, $B=B_1 \otimes B_1$.
Then $\eta _{22}^1 = \left[ {\begin{array}{*{20}{c}}
0&0&0&0&1&{ - 1}&0&{ - 1}&1
\end{array}} \right]$ and $\eta _{22}^2 = \left[ {\begin{array}{*{20}{c}}
0&{ - 1}&1&1&0&{ - 1}&{ - 1}&1&0
\end{array}} \right]$. If $c_1=0$, one has $(c_1 \eta_{22}^1+c_2 \eta_{22}^2)(B_1 \otimes B_1)=c_1=0$.
It then follows from Theorem \ref{theorem7} that $(A({\mathcal G}),B)$ is uncontrollable. This example fully demonstrates the effectiveness of the condition.

\end{example}

\begin{remark}
Theorem $6$ in \cite{chapmanconference} has established a controllability condition for Kronecker product {networks}, which requires $A(\mathcal{G})$ to be diagonalizable. This assumption is conservative and the condition can only be applied to a very restricted type of networks. The new condition here removes this restriction, thus is more general and flexible. This nontrivial extension is of great significance in engineering {applications}.
\end{remark}

In the following, some more intuitive and easily-verifiable conditions are presented, {which reveal how the controllability of the factor networks affects the controllability of the whole composite network.}

\begin{corollary}
\label{corollary2}
If the composite network (\ref{eq3})-(\ref{eq4}) is controllable, then the factor networks $(A(\mathcal{G}_1),B_1)$ and $(A(\mathcal{G}_2),B_2)$ are controllable.

\end{corollary}

\begin{corollary}
\label{corollary3}
If $A(\mathcal{G}_1)$ has an eigenvalue $0$, to ensure the controllability of the composite network (\ref{eq3})-(\ref{eq4}), it is necessary that $B_2=I_n$. Moreover, if $A(\mathcal{G}_2)$ has an eigenvalue $0$, to ensure the controllability of the composite network, it is necessary that $B_1=I_N$.

\end{corollary}

{\textbf{Proof:}}
If $A(\mathcal{G}_1)$ has an eigenvalue $0$, {without loss of generality, let} $\lambda_i=0$, then the left eigenvectors of $A(\mathcal{G})$ corresponding to $0$ are $v_i^1\otimes w_1^1$, $\cdots$, $v_i^1 \otimes w_1^{q_1}$, $\cdots$, $v_i^1\otimes w_t^1$, $\cdots$, $v_i^1 \otimes w_t^{q_t}$, where $v_i^1$ is the left eigenvector of $A(\mathcal{G}_1)$ associated with the eigenvalue $0$; $w_1^1$, $\cdots$, $w_1^{q_1}$, $\cdots$, $w_t^1$, $\cdots$, $w_t^{q_t}$ are all the left root vectors of $A(\mathcal{G}_2)$. If the composite network (\ref{eq3})-(\ref{eq4}) is controllable, then $\left[ {\sum\limits_{j = 1}^t {\sum\limits_{k = 1}^{{q_k}} ({c_j^kv_i^1 \otimes w_j^k)} } } \right](B_1  \otimes B_2) \ne 0$, for any scalars $c_1^1$, $\cdots$, $c_1^{q_1}$, $\cdots$, $c_t^1$, $\cdots$, $c_t^{q_t}$, which are not all zero. That is, $(v_i^1 B_1 ) \otimes \left[ {\sum\limits_{j = 1}^t {\sum\limits_{k = 1}^{{q_k}} {(c_j^kw_j^kB_2)} } } \right] \ne 0$, which implies $\left[ {\sum\limits_{j = 1}^t {\sum\limits_{k = 1}^{{q_k}} {c_j^kw_j^k} } } \right]B_2 \ne 0$, for any {nonzero $[c_1^1,\cdots, c_1^{q_1}, \cdots, c_t^1, \cdots, c_t^{q_t}]$}. Since $w_1^1$, $\cdots$, $w_1^{q_1}$, $\cdots$, $w_t^1$, $\cdots$, $w_t^{q_t}$ are all the left root vectors of $A(\mathcal{G}_2)$, which span $\mathbb{C}^n$, one has $B_2=I_n$.

The second part can be proved similarly.
\hfill $\blacksquare$

Hereinafter, the special case where $A(\mathcal{G}_1)$ or $A(\mathcal{G}_2)$ is diagonalizable is analyzed in detail, for which a simple and effective condition is obtained. Compared with the conditions for the general case, this easy-to-verify condition {allows} to check the {network} controllability {more efficiently}.

\subsection{$A(\mathcal{G}_1)$ or $A(\mathcal{G}_2)$ is diagonalizable}

In this subsection, a controllability criterion is established for the Kronecker product network {(\ref{eq3})-(\ref{eq4})} with one diagonalizable factor network.

\begin{corollary}
\label{theorem11}
Assume that $A(\mathcal{G}_1)$ or $A(\mathcal{G}_2)$ is diagonalizable. Let $\Lambda{=}\{\lambda_1,{\cdots},\;\lambda_s\}$ be the set of the eigenvalues of $A(\mathcal{G}_1)$, and $\Sigma=\{\mu_1,\cdots,\mu_t\}$ be the set of the eigenvalues of $A(\mathcal{G}_2)$. The composite network (\ref{eq3})-(\ref{eq4}) is controllable if and only if the following four conditions hold simultaneously:
\begin{enumerate}
  \item $(A(\mathcal{G}_1),B_1)$ and $(A(\mathcal{G}_2),B_2)$ are controllable;
  \item if $A(\mathcal{G}_1)$ has an eigenvalue $0$, then $B_2=I_n$;
  \item if $A(\mathcal{G}_2)$ has an eigenvalue $0$, then $B_1=I_N$;
  \item if $\lambda_{i_1}\mu_{j_1}=\lambda_{i_2}\mu_{j_2}=\cdots=\lambda_{i_r}\mu_{j_r} \ne 0$, where $\lambda_{i_k} \in \Lambda$, $\mu_{j_k} \in \Sigma$, for $k=1,\cdots,r$, $r >1$, then $({v_{{i_1}}}B_1 ) \otimes (w_{{j_1}}B_2)$, $\cdots$, $({v_{{i_r}}}B_1 ) \otimes ({w_{{j_r}}}B_2)$ are linearly independent, where $v_{i_k}$ is the left eigenvector of $A(\mathcal{G}_1)$ corresponding to the eigenvalue $\lambda_{i_k}$; $w_{j_k}$ is the left eigenvector of $A(\mathcal{G}_2)$ corresponding to the eigenvalue $\mu_{i_k}$, $k=1,\cdots,r$.
\end{enumerate}

\end{corollary}

{\textbf{Proof:}}
The case that $A(\mathcal{G}_1)$ is diagonalizable is firstly proved. In this case, $s=N$. Then with a similar method, one can prove the other case easily.

%
%
%
%
%
%
Necessity: First of all, from Corollaries \ref{corollary2} and \ref{corollary3}, it follows that conditions $(1)$-$(3)$ are necessary for the controllability of the composite network (\ref{eq3})-(\ref{eq4}).

Assume that $\lambda_{i_1}\mu_{j_1}=\lambda_{i_2}\mu_{j_2}=\cdots=\lambda_{i_r}\mu_{j_r}=\sigma \ne 0$, where $\lambda_{i_k} \in \Lambda$, $\mu_{j_k} \in \Sigma$, $k=1,\cdots,r$, $r >1$. Then, {any} left eigenvector of $A(\mathcal{G})$ corresponding to $\sigma$ can be expressed in the form of $\sum\limits_{k = 1}^r {{\alpha _k}({v_{{i_k}}} \otimes {w_{{j_k}}})} $, where $v_{i_k}$ is the left eigenvector of $A(\mathcal{G}_1)$ corresponding to the eigenvalue $\lambda_{i_k}$; $w_{j_k}$ is the left eigenvector of $A(\mathcal{G}_2)$ corresponding to the eigenvalue $\mu_{i_k}$; $\alpha_{k} \in \mathbb{R}$ $(k=1,\cdots,r)$ are {some scalars} which are not all zero. If the composite network (\ref{eq3})-(\ref{eq4}) is controllable, then $\left[ {\sum\limits_{k = 1}^r {{\alpha _k}({v_{{i_k}}} \otimes {w_{{j_k}}})} } \right](B_1  \otimes B_2) \ne 0$. Consequently,
${ {\sum\limits_{k = 1}^{{r}} {\alpha _{k}({v_{i_k}}{B_1}) \otimes (w_{j_k}{B_2})} } }  \ne 0$,
{for any scalars $\alpha _1$, $\cdots$, $\alpha_{r}$, which are not all zero.}
Therefore, $({v_{{i_1}}}B_1 ) \otimes ({w_{{j_1}}}B_2)$, $\cdots$, $({v_{{i_r}}}B_1 ) \otimes ({w_{{j_r}}}B_2)$ are linearly independent.


Sufficiency: One needs to prove that, if the composite network (\ref{eq3})-(\ref{eq4}) is uncontrollable, then at least one condition in Corollary \ref{theorem11} does not hold. Based on {the PBH test}, if the composite network (\ref{eq3})-(\ref{eq4}) is uncontrollable, then there exists a left-eigenpair of $A(\mathcal{G})$, denoted as $(\sigma, \xi )$, such that
\begin{equation}
\label{eq90}
\xi B=0.
\end{equation}
This will be discussed in the following $3$ cases.
\begin{itemize}
  \item If {$\sigma \ne 0$ is an eigenvalue with the geometric multiplicity being $1$, assume that $\sigma=\lambda_{i}\mu_{j}$, it then follows from Theorem \ref{theorem6} that $\xi=v_{i}\otimes w_{j}^1$,} where $\lambda_i \in \Lambda\ne 0 $ with the corresponding left eigenvector $v_i$; $\mu_j  \in \Sigma \ne 0$ with the corresponding left eigenvector $w_j^1$. From {equality} (\ref{eq90}), it can be easily {deduced} that $(v_{i}\otimes w_{j}^1)(B_1 \otimes B_2)=(v_i B_1)\otimes (w_j^1 B_2) =0$, which yields
$v_i B_1=0$
or
$w_j^1 B_2=0$.
It can be then derived that $(A(\mathcal{G}_1),B_1)$ is uncontrollable or $(A(\mathcal{G}_2),B_2)$ is uncontrollable.
  \item If $\sigma \ne 0$ is {an eigenvalue with the geometric multiplicity more than $1$}, assume that $\lambda_{i_1}\mu_{j_1}=\cdots=\lambda_{i_r}\mu_{j_r}=\sigma$, where $\lambda_{i_k} \in \Lambda \ne 0$, $\mu_{j_k}  \in \Sigma \ne 0$, $k=1,\cdots,r$, $r >1$. It follows from Theorem \ref{theorem6} that $\xi$ can be expressed in the form of $\sum\limits_{k = 1}^r {{\hat \alpha _k}({v_{{i_k}}} \otimes {w_{{j_k}}})} $, where ${v_{{i_k}}}$ is the left eigenvector of $A(\mathcal{G}_1)$ associated with the eigenvalue ${{\lambda}_{i_k}}$; ${w_{{j_k}}}$ is the left eigenvector of $A(\mathcal{G}_2)$ associated with the eigenvalue ${{\mu}_{j_k}}$; $\hat \alpha_{k} \in \mathbb{R}$ $(k=1,\cdots,r)$ are some scalars, which are not all zero. {From equality (\ref{eq90}), it follows that} there exists a nonzero vector $[\hat \alpha_1,\cdots,\hat \alpha_r]$, such that
      $\left[ {\sum\limits_{k = 1}^r {{{\hat \alpha }_k}({v_{{i_k}}} \otimes {w_{{j_k}}})} } \right](B_1  \otimes B_2) = \sum\limits_{k = 1}^r {{{\hat \alpha }_k}({v_{{i_k}}}B_1 ) \otimes ({w_{{j_k}}}\\B_2)}  = 0.$
      It indicates that $({v_{{i_1}}}B_1 ) \otimes ({w_{{j_1}}}B_2)$, $\cdots$, $({v_{{i_r}}}B_1 ) \otimes ({w_{{j_r}}}B_2)$ are linearly dependent.

  \item If $\sigma = 0$, assume that $\lambda_{1}=\cdots=\lambda_{p}=0$ and $\mu_1=\cdots=\mu_l=0$. It follows from Theorem \ref{theorem6} that $\xi$ can be expressed in the form of $\sum\limits_{i = 1}^p {\sum\limits_{j = 1}^t {\sum\limits_{k = 1}^{{q_j}} ({\hat \alpha _{ij}^k{v_i} \otimes w_j^k}) } }  + \sum\limits_{j = 1}^l {\sum\limits_{i = p + 1}^N {(\hat \beta_{ij}{v_i} \otimes w_j^1} }) $, where ${v_{i}}$ is the left eigenvector of $A({\mathcal{G}_1})$ associated with the eigenvalue $\lambda_i$, {$i=1,\cdots,N$}; $q_j$, ${w_{{j}}^k}$ ($k=1,\cdots,q_j$) are defined as in Theorem \ref{theorem6}; ${\hat \alpha _{ij}^k}$ $(k=1,\cdots,q_j,\;j=1,\cdots,t,\;i=1,\cdots,p)$ and ${\hat \beta _{ij}}$ $(j=1,\cdots,l,\;i=p+1,\cdots,N)$ are some scalars, which are not all zero. {From equality (\ref{eq90}), it follows} that there exist scalars ${\hat \alpha _{ij}^k}$ $(k=1,\cdots,q_j,\;j=1,\cdots,t,\;i=1,\cdots,p)$ and ${\hat \beta _{ij}}$ $(j=1,\cdots,l,\;i=p+1,\cdots,N)$, which are not all zero, such that
      $\left[ {\sum\limits_{i = 1}^p {\sum\limits_{j = 1}^t {\sum\limits_{k = 1}^{{q_j}} {(\hat \alpha _{ij}^k{v_i} \otimes w_j^k)} } }  + \sum\limits_{j = 1}^l {\sum\limits_{i = p + 1}^N {(\hat \beta _{ij}{v_i} \otimes w_j^1)} } } \right]({B_1}\\ \otimes {B_2}) = 0$,
      which implies that
 \begin{equation}
      \label{eq140}
\begin{array}{l}
\sum\limits_{j = 1}^l {\left[ {(\sum\limits_{i = 1}^p {\hat \alpha _{ij}^1{v_i}{B_1}}  + \sum\limits_{i = p + 1}^N {{{\hat \beta }_{ij}}{v_i}{B_1}} ) \otimes (w_j^1{B_2})} \right]}  + \\
\sum\limits_{j = 1}^l {\sum\limits_{k = 2}^{{q_j}} \left[ {(\sum\limits_{i = 1}^p {\hat \alpha _{ij}^k{v_i}{B_1})}    \otimes (w_j^k{B_2})} \right]} + \\
\sum\limits_{j = l + 1}^t {\sum\limits_{k = 1}^{{q_j}} {\left[ {(\sum\limits_{i = 1}^p {\hat \alpha _{ij}^k{v_i}{B_1}) \otimes (w_j^k{B_2})} } \right]} }  = 0.
\end{array}
\end{equation}
      From Lemma \ref{lemma2}, it can be derived that equality (\ref{eq140}) holds if and only if at least one of the following conditions holds:
      \begin{enumerate}
        \item $w_1^1B_2$, $\cdots$, $w_1^{q_1}B_2$, $\cdots$, $w_t^1B_2$, $\cdots$, $w_t^{q_t}B_2$ are linearly dependent.
        \item  $v_1B_1$, $\cdots$, $v_NB_1$ are linearly dependent.
      \end{enumerate}
      Since $w_1^1$, $\cdots$, $w_1^{q_1}$, $\cdots$, $w_t^1$, $\cdots$, $w_t^{q_t}$ are all the left root vectors of $A(\mathcal{G}_2)$, which span $\mathbb{C}^n$, the first condition implies that $B_2 \ne I_n$. Moreover, since $v_1$, $\cdots$, $v_N$ are all the left eigenvectors of $A(\mathcal{G}_1)$, which span $\mathbb{C}^{N}$, the second condition indicates that $B_1 \ne I_N$. Thus, in this case, if the composite network (\ref{eq3})-(\ref{eq4}) is uncontrollable, then $B_1 \ne I_N$ or $B_2 \ne I_n$.
\end{itemize}

Therefore, if the composite network (\ref{eq3})-(\ref{eq4}) is uncontrollable, then at least one condition in Corollary \ref{theorem11} does not hold.
This completes the proof.
\hfill $\blacksquare$

The effectiveness of this criterion can be illustrated by the following example.

\begin{example}
Consider two graphs ${\mathcal G}_1$ and ${\mathcal G}_2$, which are depicted in Fig.\ref{fig2}.
\begin{figure}[htb]
\centering
\subfloat[${\mathcal G}_1$]{%
\includegraphics[width=.2\textwidth]{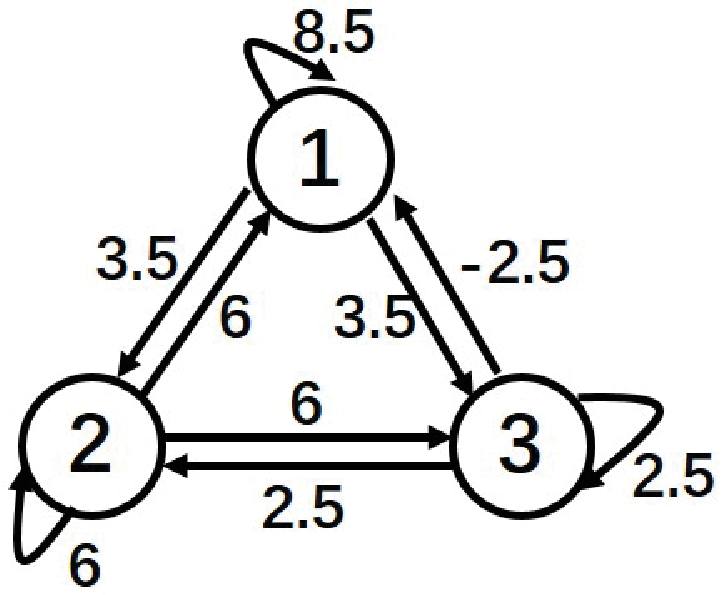}}
\subfloat[${\mathcal G}_2$]{%
\includegraphics[width=.2\textwidth]{G_ex.eps}}\hfill
\caption{Factor graphs ${\mathcal G}_1$ and ${\mathcal G}_2$}
\label{fig2}
\end{figure}
It is easy to verify that
$A({{\mathcal G}_1}) = \left[ {\begin{array}{*{20}{c}}
   8.5 & 6 & -2.5  \\
   3.5 & 6 & 2.5  \\
   3.5 & 6 & 2.5  \\
\end{array}} \right]$ and $A({{\mathcal G}_2}) = \left[ {\begin{array}{*{20}{c}}
   8.5 & 4 & -0.5  \\
   3.5 & 8 & 0.5  \\
   3.5 & 3 & 5.5  \\
\end{array}} \right]$.
Note that $A({\mathcal G}_1)$ is diagonalizable. The eigenvalues of $A({\mathcal G}_1)$ are $\lambda_1=12$, $\lambda_2=5$ and $\lambda_3=0$ with the corresponding left eigenvectors $v_1=e_1+e_2$, $v_2=-e_1+e_3$ and $v_3=-e_2+e_3$, respectively. Moreover, the eigenvalues of $A({\mathcal G}_2)$ are $\mu_1=12$, $\mu_2=5$. The left eigenvector associated with $\mu_1$ is $w_1=-e_1-e_2$ and the left Jordan chain corresponding to $\mu_2$ is $w_2^1=e_2-e_3$, $w_2^2=-e_1+e_2$.
Assume that the first two nodes of ${{\mathcal G}_1}$ have control inputs and all the nodes of ${{\mathcal G}_2}$ are under control.
Thus, one has
${B_1} = [e_1^T,\;e_2^T]$ and ${B_2} = [e_1^T,\;e_2^T,\;e_3^T]$. Firstly, it is easy to verify that $(A({\mathcal G}_1),B_1)$ and $(A({\mathcal G}_2),B_2)$ are controllable. Since $A({\mathcal G}_1)$ has an eigenvalue $0$ and $B_2=I_3$, the second condition holds. Further, noting that $\lambda_1\mu_2=\lambda_2\mu_1=60$, one needs to check whether $(v_1B_1)\otimes (w_2^1B_2)$ and $(v_2B_1)\otimes (w_1B_2)$ are linearly dependent. Since $B_2=I_3$, $w_2^1B_2$ and $w_1B_2$ are linearly independent. Moreover, $v_1B_1 \ne 0$ and $v_2B_1 \ne 0$. Thus, one can verify that $(v_1B_1)\otimes (w_2^1B_2)$ and $(v_2B_1)\otimes (w_1B_2)$ are linearly independent. It then follows from Corollary \ref{theorem11} that $(A({\mathcal G}),B)$ is controllable. This result coincides with the conclusion derived by the Kalman rank condition.

\end{example}

If $A(\mathcal{G})$ is diagonalizable, then both $A(\mathcal{G}_1)$ and $A(\mathcal{G}_2)$ are diagonalizable. The {conditions} in Corollary \ref{theorem11} {are} also effective for diagonalizable Kronecker product {networks}, {as demonstrated by the following example}.

\begin{example}
\label{example1}
Consider two graphs ${\mathcal G}_1$ and ${\mathcal G}_2$, which are depicted in Fig. \ref{ex_fig}.
\begin{figure}[htb]
\centering
\subfloat[${\mathcal G}_1$]{%
\includegraphics[width=.2\textwidth]{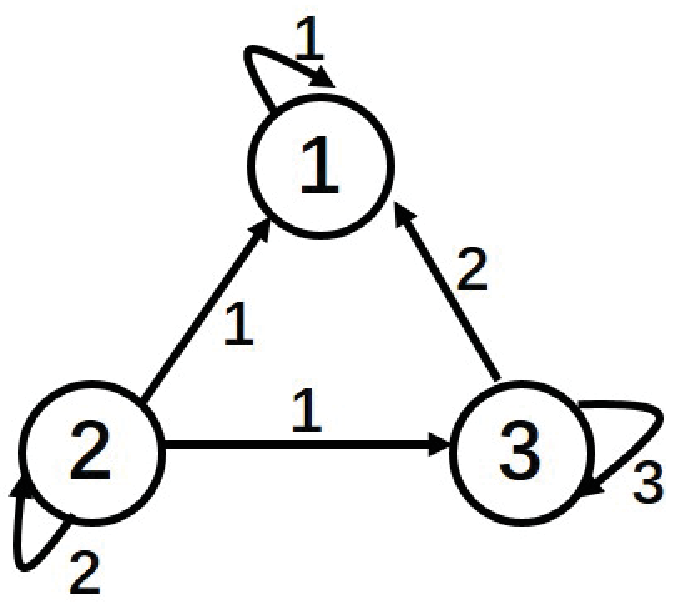}}
\subfloat[${\mathcal G}_2$]{%
\includegraphics[width=.2\textwidth]{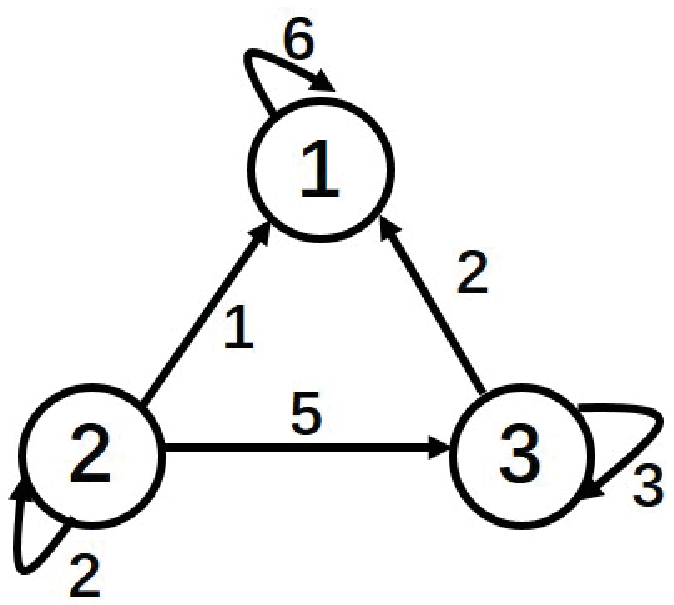}}\hfill
\caption{Factor graphs ${\mathcal G}_1$ and ${\mathcal G}_2$}
\label{ex_fig}
\end{figure}
It is easy to verify that
$A({{\mathcal G}_1}) = \left[ {\begin{array}{*{20}{c}}
   1 & 1 & 2  \\
   0 & 2 & 0  \\
   0 & 1 & 3  \\
\end{array}} \right]$ and $A({{\mathcal G}_2}) = \left[ {\begin{array}{*{20}{c}}
   6 & 1 & 2  \\
   0 & 2 & 0  \\
   0 & 5 & 3  \\
\end{array}} \right]$.
The eigenvalues of $A({\mathcal G}_1)$ are $\lambda_1=2$, $\lambda_2=1$ and $\lambda_3=3$ with the corresponding left eigenvectors
$v_1=e_2$, $v_2=e_1-e_3$ and $v_3=e_2+e_3$, respectively.
The eigenvalues of $A({\mathcal G}_2)$ are $\mu_1=2$, $\mu_2=3$ and $\mu_3=6$ with the corresponding left eigenvectors $w_1=e_2$, $w_2=5e_2+e_3$, $w_3=3/2e_1+13/8e_2+e_3$, respectively. It is easy to verify that $A({\mathcal G})$ is diagonalizable.
Assume that the first two nodes of ${{\mathcal G}_1}$ have control inputs. Node $2$ and node $3$ of ${{\mathcal G}_2}$ are under control.
Thus, one has
${B_1} = [e_1^T,\;e_2^T]$ and ${B_2} = [e_2^T,\;e_3^T]$.
Both pairs $(A({\mathcal G}_1),B_1)$ and $(A({\mathcal G}_2),B_2)$ are controllable. Since $\lambda_1\mu_2=\lambda_2\mu_3=\lambda_3\mu_1=6$,
the value of $[c_1(v_1 \otimes w_2)+c_2(v_2\otimes w_3)+c_3(v_3 \otimes w_1)](B_1 \otimes B_2)$ is required to be checked, which gives $[{c_1}({v_1} \otimes {w_2}) + {c_2}({v_2} \otimes {w_3}) + {c_3}({v_3} \otimes {w_1})]({B_1} \otimes {B_2}){\rm{ = }}[ {\begin{array}{*{20}{c}}
{\frac{{{\rm{13}}}}{{\rm{8}}}{c_2}},&{{c_2}},&{5{c_1} + {c_3}},&{{c_1}}
\end{array}} ] \ne 0$, for any nonzero $[c_1,\;c_2,\;c_3]$.
It then follows from Corollary \ref{theorem11} that $(A({\mathcal G}),B)$ is controllable, which coincides with the conclusion derived by the Kalman rank condition. This example fully demonstrates the effectiveness of the condition.

\end{example}

\begin{remark}
A controllability criterion for diagonalizable Kronecker product {networks} was established in \cite{chapmanconference},
which was claimed to be necessary and sufficient. However, Example \ref{example1} shows that the necessity of that criterion does not hold. Corollary \ref{theorem11}
provides a new controllability condition, which is not only sufficient but also necessary.

\end{remark}
%
%
%
%
%
%
%
%
%
%
%
%
%
%
%
%
%
%
%
%
\section{Controllability of Higher-dimensional Multi-agent Systems Revisited}
\label{section5}

\subsection{Problem Statement}
Consider a multi-agent system consisting of $N$ agents labeled by the set $V = \left\{ {1,2, \cdots ,N} \right\}$. Assign the roles of leaders and followers to the agents by denoting ${V_L} = \left\{ {{v_1},{v_2}, \cdots ,{v_m}} \right\}$ and ${V_F} = V\backslash {V_L}$ {as} the sets of indices of the leaders and followers, respectively, where $m \le N$.

To each follower $i \in V_F$, associate a dynamical system
${{\dot x}_i} = {z_i}$,
and to each leader $i \in V_L$, we associate a dynamical system
${{\dot x}_i} = {z_i}{\rm{ + }}B{u_i}$,
where ${x_i} \in {\mathbb{R}^n}$ is the state of agent $i$; ${u_i} \in {\mathbb{R}^p}$ is the external input to agent $i \in V_L$; $z_i \in {\mathbb{R}^n}$ is the coupling input from other agents.

Agent $i$ is said to be a neighbor of agent $j$ if its state is known by agent $j$. Here, assume that the neighboring relationships are fixed, {which} can be described by a directed and weighted graph ${\mathcal{G}}=(V,E,W)$. The coupling input $z_i$ to each agent $i \in V$ is determined by the diffusive coupling rule based on the neighboring relations as follows:
${z_i} = H\sum\limits_{(j,i) \in E} {w_{ij}({x_j} - {x_i})}$,
where $H \in \mathbb{R}^{n \times n}$ is the matrix describing inner-coupling between different components, and $w_{ij} \in \mathbb{R}$ is the strength of the information link.


The Laplacian matrix of $\mathcal{G}$ and the external input channels of the multi-agent system are denoted by
$L  \in {\mathbb{R}^{N \times N}}$ and $\Delta = diag\left\{{d_1},{d_2}, \cdots, {d_N}\right\}$,
respectively, where {${d_i} = 1$} for $i \in V_L$, but otherwise {${d_i}=0$}, for all $i = 1,2, \cdots ,N$.
Let $X = {[ {x_{\rm{1}}^T,\;x_2^T,\; \cdots ,\;x_N^T} ]^T}$ be the whole state of the {multi-agent} system, and $U = {[ {u_{\rm{1}}^T,\;u_2^T,\; \cdots ,\;u_N^T} ]^T}$ be the total external control input. Then, the above multi-agent system can be rewritten in a compact form as
{\begin{equation}
\label{eq14}
 \begin{array}{l}
 \dot X = F X + G U,
 \end{array}
\end{equation}}
with
{\begin{equation}
\label{eq15}
F  =  -L \otimes H,\;\;G  = \Delta \otimes B.
\end{equation}}
Note that matrix $F$ {has} the form of the Kronecker product of two matrices and so {does} $G$. This higher-dimensional multi-agent system can be seen as a special case with $A(\mathcal{G}_1)$ being a Laplacian matrix. In the following,
conditions for ensuring the controllability of the multi-agent system (\ref{eq14})-(\ref{eq15}) are specified.
{\begin{remark}
The controllability of networked LTI systems or multi-agent systems with linear dynamics has been investigated in \cite{caoming_2014,hao_controllability_2017,xuemengran_comment_2018,xuemengran_input_2017,xuemengran_model_2018}. The state matrices for those systems have the form of $I \otimes A+ L \otimes H$ rather than a pure Kronecker product. The higher-dimensional multi-agent systems investigated here can be seen as a special case with nodes having no internal dynamics.
\end{remark}}
\subsection{A Counterexample}
Recall the controllability condition {for} the multi-agent system (\ref{eq14})-(\ref{eq15}) given in \cite{caining}, where it was assumed that the first $N_l$ agents are leaders. Then, $\Delta  = [ {\begin{array}{*{20}{c}}
{{e_1^T}},& \cdots, &{{e_{{N_l}}^T}},&\textbf{0}_N,& \cdots, &\textbf{0}_N
\end{array}} ] \in {\mathbb{R}^{N \times N}}$ and $L$ can be partitioned as
$L = \left[ {\begin{array}{*{20}{c}}
{{L_{ll}}}&{{L_{lf}}}\\
{{L_{fl}}}&{{L_{ff}}}
\end{array}} \right]$,
where $L_{ll} \in \mathbb{R}^{N_l \times N_l}$ and $L_{ff} \in \mathbb{R}^{(N-N_l) \times (N-N_l)}$ with subscripts `l' and `f' denoting `leader' and `follower', respectively.

\begin{theorem}
\label{error}
(Theorem $1$ of \cite{caining})
For an LTI swarm system described by (\ref{eq14})-(\ref{eq15}), suppose the first $N_l$ agents are leaders. Then, the system is completely controllable if and only if
\begin{enumerate}
  \item $(H,B)$ is a controllable matrix pair;
  \item $L$ represents a controllable graph, i,e, $[ {\begin{array}{*{20}{c}}
{{L_{ff}}}&{{L_{fl}}}
\end{array}} ]$ is a controllable matrix pair.
\end{enumerate}

\end{theorem}

However, while this condition is necessary for the controllability of the network, it may not be sufficient, as shown in the following example.

\begin{figure}[!htb]
  \centering
  \includegraphics[scale=0.4]{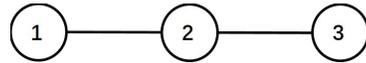}
  \caption{A path network}
  \label{fig4}
\end{figure}

Consider the undirected path network with three nodes {depicted} in Fig. \ref{fig4}.
Node $1$ is selected to be the leader. One has
$L  = \left[ {\begin{array}{*{20}{c}}
1&{ - 1}&0\\
{ - 1}&2&{ - 1}\\
0&{ - 1}&1
\end{array}} \right]$ and $\Delta  = \left[ {\begin{array}{*{20}{c}}
1&{ 0}&0\\
{ 0}&0&{ 0}\\
0&{ 0}&0
\end{array}} \right]$. Let $H = \left[ {\begin{array}{*{20}{c}}
{1.5}&{0.5}\\
{0.5}&{1.5}
\end{array}} \right]$ and $B = \left[ {\begin{array}{*{20}{c}}
1\\
2
\end{array}} \right]$. It is easy to verify that $(H,B)$ and $(L_{ff}, L_{fl})$ are controllable. From Theorem $1$ in \cite{caining}, it follows that this undirected path network is controllable.

However, if one checks the controllability of this system by using the classical Kalman rank condition, one can find that $(-L \otimes H, \Delta \otimes B)$ is actually uncontrollable. Therefore, its sufficiency does not hold.

{This example is a modification of the example presented in \cite{xuemengran_comment_2018}, which was used to demonstrate that the controllability condition for diffusive networks proposed in \cite{caoming_2014} is not always sufficient. The {incomplete} eigenanalysis of networks presented in \cite{caoming_2014} and \cite{caining} leads to errors in controllability analysis. In the following, a modified controllability condition is proposed, which is both necessary and sufficient.}

\subsection{New Controllability Criteria}
{Note that the conditions proposed in Section \ref{section4} are not restricted to checking the controllability of Kronecker product networks. They can tackle the controllability problem for any system represented by the Kronecker product of two matrices. Based on the results in Section \ref{section4},} a modified controllability condition for the multi-agent system (\ref{eq14})-(\ref{eq15}) is established as follows.

\begin{corollary}
\label{corollary12}
The multi-agent system (\ref{eq14})-(\ref{eq15}) is controllable if and only if the following three conditions hold simultaneously:
\begin{enumerate}
  \item $(L,\Delta)$ is controllable;
  \item $rank(B)=n$;
  \item if $H$ has an eigenvalue $0$, then $\Delta=I_N$.
\end{enumerate}

\end{corollary}
{\textbf{Proof:}
Note that $L$ has an eigenvalue $0$. If the multi-agent system (\ref{eq14})-(\ref{eq15}) is controllable, then $rank(B)=n$.
From Corollaries \ref{corollary2} and \ref{corollary3}, it follows that the conditions (1) and (3) are necessary for the controllability of the multi-agent system (\ref{eq14})-(\ref{eq15}).}

{For sufficiency, one needs to prove that, if the conditions (1)-(3) hold, then no left eigenvectors of $F$ are orthogonal to $G$, thus the multi-agent system (\ref{eq14})-(\ref{eq15}) is controllable. It is easy to verify that if the conditions (1)-(3) hold, then no left eigenvectors of $F$ associated with eigenvalue $0$ are orthogonal to $G$. For a nonzero eigenvalue $\sigma=-\lambda\mu$, if it is not a common eigenvalue, it follows from Theorem \ref{theorem6} that any corresponding left eigenvector can be expressed as $\xi=\left(\sum\limits_{i = 1}^{{\theta}} {{c_i}v^i}\right)  \otimes w^1 + \sum\limits_{k = 2}^{{\theta }} {\left[ {\left( {\sum\limits_{i = k}^{{\theta}} {\sum\limits_{j = 1}^{i - k + 1} {{c_i}{l_{ij}}v^j} } } \right) \otimes w^k} \right]}$, where $v^{1}$, $\cdots$, $v^{p}$ is the left Jordan chain of $L$ associated with the eigenvalue $\lambda$, and $w^{1}$, $\cdots$, $w^{q}$ is the left Jordan chain of $H$ associated with the eigenvalue $\mu$, $\theta=\min\{p,q\}$, $l_{ij}$ is a nonzero scalar about $\lambda$ and $\mu$, $c_1$, $\cdots$, $c_{\theta}$ are scalars, which are not all zero. Next, the value of $\xi(\Delta \otimes B)$ will be checked.
It is easy to verify that
$\xi (\Delta  \otimes B) = \left(\sum\limits_{i = 1}^{{\theta}} {{c_i}v^i\Delta }\right) \otimes (w^1B) + \sum\limits_{k = 2}^{{\theta }} {\left[ {\left( {\sum\limits_{i = k}^{{\theta }} {\sum\limits_{j = 1}^{i - k + 1} {{c_i}{l_{ij}}v^j\Delta } } } \right) \otimes (w^kB)} \right]}$.
Since $rank(B)=n$, $w^1B$, $\cdots$, $w^{\theta}B$ are linearly independent. Since $(L,\Delta)$ is controllable, it follows that $\sum\limits_{i = 1}^{{\theta}} {{c_i}v^i\Delta }$ and ${\sum\limits_{i = k}^{{\theta }} {\sum\limits_{j = 1}^{i - k + 1} {{c_i}{l_{ij}}v^j\Delta } } } $ $(k=2,\cdots,\theta)$ are not all zero. From Lemma \ref{lemma2}, it can be easily deduced that $\xi (\Delta  \otimes B)\ne 0$. Thus, for any nonzero eigenvalue of $F$, if it is not a common eigenvalue, no corresponding left eigenvectors are orthogonal to $G$.
Similarly, one can prove that, for each nonzero common eigenvalue of $F$, no corresponding left eigenvectors are orthogonal to $G$.
Consequently, no left eigenvectors of $F$ are orthogonal to $G$.
Therefore, the multi-agent system (\ref{eq14})-(\ref{eq15}) is controllable. This completes the proof.}
{\hfill $\blacksquare$}

According to Corollary \ref{corollary12}, one can easily verify that the undirected path network in the above counterexample is uncontrollable. This example fully demonstrates the effectiveness of the condition.
%
%
%
%
%
%
%
%
%
%
%
%
%
%
%

\section{Conclusions}
\label{section6}
The controllability of {Kronecker product networks} has been investigated, in which the factor networks have general directed topologies. A necessary and sufficient condition for the controllability of the composite network has been derived, which is effective and has a much lower computational cost as compared to existing criteria. For the special case where at least one factor network is diagonalizable, a specified condition has also been established, which is simple and easier to verify. Moreover, the controllability of higher-dimensional multi-agent systems has been reinvestigated. It is found that the sufficiency of the controllability criterion given in \cite{caining} does not hold. Consequently, a modified condition is derived, which is necessary and sufficient. In future studies, the controllability and observability of other types of network-of-networks will be further considered.





\end{document}